\documentclass{article}
\usepackage{amsmath,amssymb,latexsym,amsthm,amscd}
\font\goth=eusm10
\newcommand\F{\mathcal F}
\newcommand\E{\mathcal E}
\newcommand\Ii{\hbox{\goth I}}

\newcommand\G{\mathcal G}

\newcommand\Oc{\hbox{\goth O}}

\newcommand\CC{\mathbf C}

\newcommand\Pj{\mathbf{P}}

\newcommand\ZZ{\mathbf{Z}}

\newtheorem{theorem}{Theorem}[section]
\newtheorem{RR}[theorem]{Theorem}
\newtheorem{proposition}[theorem]{Proposition}
\newtheorem{definition}[theorem]{Definition}

\newtheorem{question}[theorem]{Question}

\newtheorem{remark}[theorem]{Remark}

\newtheorem{lemma}[theorem]{Lemma}
\newtheorem*{Lemma}{Lemma}

\theoremstyle{plain}
\theoremstyle{definition}
\theoremstyle{remark}
\numberwithin{equation}{section}
\begin{document}
\title
{Rank 4 vector bundles on the
quintic threefold}
\author{Carlo Madonna
\footnote{Dipartimento di Matematica,
Universit\`a degli Studi di Roma "La Sapienza",
P.le A.Moro 1, 00185 Roma, Italia. email: madonna@mat.uniroma1.it}}
\date{}
\maketitle
\begin{abstract}
By the results of the author and Chiantini in \cite{cm1},
on a general quintic threefold
$X \subset \Pj^4$ the minimum integer $p$ for which
there exists a positive dimensional family of irreducible rank
$p$ vector bundles on $X$ without intermediate cohomology
is at least three.
In this paper we show that $p \leq 4$,
by constructing series of positive dimensional families
of rank $4$ vector bundles on $X$ without intermediate cohomology.
The general member of such family is an indecomposable
bundle from the extension class $Ext^1(E,F)$, for a suitable
choice of the rank $2$ ACM bundles $E$ and $F$ on $X$.
The existence of such bundles of rank $p = 3$ remains under question.
\end{abstract}

\section{Introduction} \label{S:intro}

Let $X \subset \Pj^4$ be a smooth quintic hypersurface
and let $E$ be a rank $2$ vector bundle without intermediate
cohomology, i.e. such that
\begin{equation} \label{eq:int}
h^i(X,E(n))=0
\end{equation}
for all $n \in \ZZ$ and $i=1,2$.
In \cite{mad3} we found all the possible Chern classes of
an indecomposable rank $2$ vector bundle
satisfying condition (\ref{eq:int}). Moreover in \cite{cm1}
we showed, when $X$ is general, if such bundles exist
then they are all infinitesimally rigid, i.e. $Ext^1(E,E)=0$.

On the other hand
it was showed in \cite{BGS} the existence of infinitely many isomorphism
classes of irreducible vector bundles without intermediate
cohomology on any smooth hypersurface $X_r$
of degree $r \geq 3$ in $\Pj^4$.
It can be checked that when the hypersurface is general then
the rank of these bundles
is $2^3$.
Hence we introduced in \cite{cm1} the number
$$
BGS(X_r)
$$
defined as the minimum positive integer $p$ for which
there exists a positive dimensional family
of irreducible rank $p$ vector bundles without intermediate cohomology on
$X_r$.

Then combining the above quoted results we get, on a general
quintic $X$, that
\begin{equation}
3 \leq BGS(X) \leq 8.
\end{equation}

In this paper we show the following:

\begin{theorem} \label{th:main}
If $X$ is general then $BGS(X) \leq 4$.
\end{theorem}

We should then answer the following:

\begin{question}
Let $X$ be a general quintic hypersurface in $\Pj^4$.
Could it be $BGS(X)=3$?
\end{question}

To show our main result
we give examples of rank $4$ vector bundles without intermediate
cohomology, which are not infinitesimally rigid.

The examples are constructed by means of extension
classes
\begin{equation}
0 \to E_1 \to \E \to E_2 \to 0
\end{equation}
i.e. elements in $Ext^1(E_2,E_1)$,
where $E_1$ and $E_2$ are rank $2$ bundles on $X$.
When the bundles $E_1$ and $E_2$ are not split then
$\E$ has not trivial summand.
Moreover for a suitable choice of bundles $E_1$ and $E_2$, there
exists a non trivial extension class such that
the rank $4$ bundle $\E$ which corresponds to this
class does not split as a direct
sum of two rank $2$ bundles, for reason of Chern classes.
Of course if $E_1$ and $E_2$ have no intermediate cohomology
it is so also for $\E$.
We then conclude by direct calculations to
make the right choice of bundles $E_1$ and $E_2$.

\section{Generalities}

We work over the complex numbers $\CC$ and we denote by $X \subset \Pj^4$
a smooth hypersurface of degree $5$ in $\Pj^4$.
Since $Pic(X) \cong \ZZ[H]$ is generated by the class of a hyperplane
section, given the vector bundle $E$ we identify $c_1(E)$
with the integer number $c_1$ which corresponds to $c_1(E)$
under the above isomorphism. We identify $c_2$ with $\deg c_2(E)=c_2(E)\cdot H$.
If $E$ is a rank $k$ vector bundle on $X$ we denote by
$E(n)=E \otimes \Oc_X(n)$.

\begin{definition}
A rank $k$ vector bundle $E$ is called arithmetically
Cohen-Macaulay (ACM for short) if $E$ has no intermediate cohomology,
i.e.
\begin{equation}
h^i(E(n))=0
\end{equation}
for all $i=1,2$, and $n \in \ZZ$.
\end{definition}

Theorem \ref{th:main} will follow by:

\begin{proposition} \label{prop:main}
Let $X$ be a smooth quintic hypersurface in $\Pj^4$.
Then, there exist indecomposable
rank $2$ vector bundles $E_1$
and $E_2$ on $X$ without intermediate cohomology
such that
there exists an open subset
of a positive dimensional projective space parameterizing extension
classes $Ext^1(E_2,E_1)$ which correspond to infinitely many
isomorphism classes of irreducible rank $4$ vector bundles $\E$ on $X$
without intermediate cohomology.
\end{proposition}

A proof of previous proposition will be given in the next section.

We will frequently use the following version of Riemann-Roch theorem for vector bundles:

\begin{RR} \label{thm:rr}
If $\E$ is a rank $2$ vector bundle
on a smooth hypersurface $X \subset \Pj^4$ of degree $5$
with Chern classes $c_i(\E)=c_i \in \ZZ$ for $i=1,2$,
then
\begin{equation} \label{eq:RRk}
\begin{aligned}
\chi(\E) =\frac56 c_1^3-\frac12 c_1c_2+\frac{25}{6}c_1
\end{aligned}
\end{equation}
\end{RR}

\section{The examples}

In this section we will give a proof of Proposition
\ref{prop:main} which is a direct consequence of Proposition
\ref{prop:list}
and Theorem \ref{thm:main} below.
As in \cite{cm1}
given a rank $2$ vector bundle $E$ we introduce the non negative integer
\begin{equation} \label{eq:norm}
b(E)=\max \{ n \mid h^0(E(-n)) \ne 0 \}.
\end{equation}
We say that the vector bundle $E$ is {\it normalized} if $b(E)=0$.
Notice that changing $E$ by $E(-b)$ we may always assume that $E$
is normalized. The rank two bundle $E$ is {\it semistable} if $2b-c_1 \leq 0$.
If $2b-c_1<0$ then $E$ is {\it stable}.

All the possible Chern classes of irreducible 
rank 2 ACM bundles are
listed in
the following (see \cite{mad3} and \cite{cm1}):

\begin{proposition} \label{prop:list}
Let $E$ be a normalized and indecomposable rank $2$ ACM bundle on a smooth
quintic $X$.
Then
$$
(c_1,c_2) \in A \cup B
$$
where
$$
A=\{ (-2,1), (-1,2), (0,3), (0,4), (0,5),
(1,4), (1,6), (1,8), (4,30) \}
$$
and $B=\{(2,\alpha), (3,20) \}$ with $\alpha=11,12,13,14$.
When $X$ is general, all the case in $A$ arise on $X$ and
moreover for all the pairs $(c_1,c_2)\in A \cup B$ the corresponding
rank $2$ ACM bundles are infinitesimally rigid i.e. $Ext^1(E,E)=0$.
\end{proposition}

\medskip

Below we shall construct examples of rank $4$ bundles
$\G$ as extensions of type
\begin{equation} \label{eq:ext}
0 \to F(m) \to \G \to E \to 0,
\end{equation}
where $m \le 0$, and $F$ and $E$ are indecomposable and normalized
rank 2 ACM bundles on $X$ with Chern classes as in
Proposition \ref{prop:list}.
Such nontrivial extensions $\G$ will exist
whenever the extension space $Ext^1(E,F(m))$
has positive dimension, i.e. $h^1(F(m) \otimes E^{\vee}) > 0$.
By the long exact sequence of cohomology of (\ref{eq:ext}),
any such extension $\G$ has vanishing intermediate cohomology
since $F$ and $E$ are ACM.
%Then we shall compute the dimension of the groups
%$Ext^1(E,F(m))$.

\begin{lemma}
Let $E$ and $F$ be two normalized and indecomposable rank 2 ACM bundles on
the smooth quintic $X$, and suppose that $h^0(F^{\vee}(c_1(E)-m))=0$
(hence $c_1(E)-c_1(F)-m<0$ since $F$ is normalized).
Then for any zero-locus $C \subset X$ of a global section of $E$
$$h^0(\Ii_C(c_1(E)) \otimes F^{\vee}(-m))=0.$$

Moreover, if $h^0(F^{\vee}(-m))=0$ (hence $-m-c_1(F)<0$ since $F$ is normalized)
then
$$h^3(F(m) \otimes E^{\vee})=h^0(E \otimes F^{\vee}(-m))=0.$$
\end{lemma}

\begin{proof}
From the tensored by $F^{\vee}(-m+c_1(E))$ ideal sheaf sequence
of $C \subset X$:
$$
0 \to \Ii_C(c_1(E)) \otimes F^{\vee}(-m) \to F^{\vee}(c_1(E)-m) \to
\Oc_C(c_1(E)) \otimes
F^{\vee}(-m) \to 0
$$
we get
$h^0(\Ii_C(c_1(E)) \otimes F^{\vee}(-m)) \leq h^0(F^{\vee}(c_1(E)-m))=0$.
The rank 2 bundle $E$ fits in the exact sequence
\begin{equation} \label{eq:kosz}
0 \to \Oc_X \to E \to \Ii_C(c_1(E)) \to 0;
\end{equation}
and after tensoring (\ref{eq:kosz})
by $F^{\vee}(-m)$ we get
$$
0 \to F^{\vee}(-m) \to E \otimes F^{\vee}(-m) \to \Ii_C(c_1(E)) \otimes
F^{\vee}(-m) \to 0.
$$
Therefore, since $h^0(F^{\vee}(-m))=h^0(\Ii_C(c_1(E)) \otimes F^{\vee}(-m))=0$
then
$$
h^0(E \otimes F^{\vee}(-m))=0,
$$
and by duality $h^3(F(m) \otimes
E^{\vee})=0$.
\end{proof}

\begin{remark} \label{rem:l1}
\normalfont{
Let $E$ and $F$ be in (\ref{eq:ext}),
and suppose that $\chi(F(m) \otimes E^{\vee})<0$.
Then by the above lemma,
the space of extensions (\ref{eq:ext}) will be no-empty
since
$h^1(F(m) \otimes E^{\vee})=h^0(F(m) \otimes E^{\vee})+
h^2(F(m) \otimes E^{\vee})-\chi(F(m) \otimes E^{\vee})>0$.

More generally the argument used here works whenever
$$
h^3(F(m) \otimes E^{\vee})<-\chi(F(m) \otimes E^{\vee}).
$$
In the following table we summarize the cases, which we are
interested in, depending on the
Chern classes of the bundles $E$ and $F$.
To get the value of $\chi(F(m) \otimes E^{\vee})$
we used Schubert package (see \cite{ks}), and then
by the Lemma we derived the lower bound for $d$.

\par
\medskip
\begin{tabular}{llllll}
Case & $(c_1(F),c_2(F))$ & $(c_1(E), c_2(E))$ & $\chi(F(m)
\otimes E^{\vee})$ & $m$ & $d$
\cr
(1) & (4,30) & (1,8) & $-14$ & $0$ & $>14$ \cr
(2) & (4,30) & (0,3) & $-6$ & $-1$ & $>6$ \cr
(3) & (4,30) & (0,4) & $-8$ & $-1$ & $>8$ \cr
(4) & (4,30) & (0,5) & $-10$ & $-1$ & $>10$ \cr
(5) & (1,8) & (0,3) & $-1$ & $0$ & $>1$ \cr
(6) & (1,8) & (0,4) & $-2$ & $0$ & $>2$ \cr
(7) & (1,8) & (0,5) & $-3$ & $0$ & $>3$
\end{tabular}
\par\medskip
}
\end{remark}

We are now ready to show the following:

\begin{theorem} \label{thm:main}
Let $X$ be a smooth quintic in $\Pj^4$, and let $E,F,m,d$
be as in the above table. Then in each of the cases $(1)-(7)$
there exists a d--dimensional parameter space of extensions
(\ref{eq:ext}), with a general element $\G$ an indecomposable
rank $4$ vector bundle on $X$ without intermediate cohomology.
\end{theorem}

\begin{proof}
For $F,E$ as in the above table, the dimension
$d =\dim Ext^1(E,F(m))$ is always $d>1$. Therefore
for such $F,E$ there exist nontrivial extensions
given by (\ref{eq:ext}), and let $\G$ be one of them.

Since $E$ and $F$ are ACM then by the cohomology sequence
of (\ref{eq:ext}) $\G$ is without intermediate cohomology,
and by Remark \ref{rem:l1} we need only to show that
$\G$ is indecomposable.
\par
Suppose the contrary, i.e. that $\G$ splits.
Then either
\par\smallskip
(i) $\G=\Oc_X(a) \oplus \G_1$ for $a \in \ZZ$ and $\G_1$ a rank 3 bundle
without intermediate cohomology,
or

(ii)  $\G = \G_1 \oplus \G_2$ for two rank 2 ACM bundles $\G_1$ and $\G_2$.
\par\smallskip
We show that under the conditions of the theorem both cases (i) and (ii)
are impossible.

\medskip

\noindent Let us start with case (i). In this case
the exact sequence (\ref{eq:ext}) reads as

\begin{equation} \label{eq:ext1}
\begin{CD}
0 @>{}>> F(m) @>f>> \Oc_X(a) \oplus \G_1 @>g>> E @>{}>> 0.
\end{CD}
\end{equation}

We use the following (see below for a proof)

\begin{Lemma} \label{cl:cl1}
Under the above conditions either $h^0E(-a)=0$ or
$h^0 F(-m-c_1(F)+a)=0$.
\end{Lemma}

Suppose $h^0(E(-a))=0$.
Then by the exact sequence (\ref{eq:ext1}) tensorized by $\Oc_X(-a)$ we have
$h^0(F(m-a))>0$.
Let $s$ be a non trivial global section of $F(m-a)$, then we have a map
$$
s:\Oc_X(a) \to F(m).
$$
Let $j:\Oc_X(a) \oplus \G_1 \to \Oc_X(a)$ be the projection.
Then we have the composition map
$$
\varphi:=j \circ f \circ s: \Oc_X(a) \to F(m) \to \Oc_X(a).
$$
Then $\varphi \in H^0 \Oc_X \cong \CC$ and hence it is either the identity map or the zero map.
If this map is the identity then $j \circ f$ is surjective and hence
$ker(j\circ f)\cong \Oc_X(b)$ for some $b \in \ZZ$.
Then we have exact sequence
$$
0 \to \Oc_X(b) \to F(m) \to \Oc_X(a) \to 0
$$
and $F(m)$, and hence also $F$, splits since $\dim Ext^1(\Oc_X(a),\Oc_X(b))=0$,
which is absurd.

Now suppose $\varphi$ is zero. Then $j \circ f$ is zero.
Thus the image of $F(m)$ in exact sequence
(\ref{eq:ext1}) is contained in $\G_1$.
Then the kernel of $g$ is contained in $\G_1$, being equal to the image of $f$.
Let $i:\Oc_X(a) \to \Oc_X(a) \oplus \G_1$ be the inclusion.
By the assumption the map $g \circ i:\Oc_X(a) \to E$ is
the zero map, which means that $\ker g$ is not contained in $\G_1$,
which is absurd.
\par\medskip
Suppose now that $h^0F(-m-c_1(F)+a)=0$
and consider the dual exact sequence of exact sequence (\ref{eq:ext1})
\begin{equation}
0 \to E^{\vee} \to \Oc_X(-a) \oplus \G_1^{\vee} \to F^{\vee}(-m) \to 0.
\end{equation}
Set $c=c_1(F)$ and $c'=c_1(E)$.
Since $E^{\vee}\cong E(-c')$ and $F^{\vee}(-m)\cong F(-c-m)$ the
above exact sequence reads as
$$
0 \to E(-c') \to \Oc_X(-a) \oplus \G_1^{\vee} \to F(-c-m) \to 0.
$$
This exact sequence tensorized by $\Oc_X(a)$ reads as
$$
0 \to E(-c'+a) \to \Oc_X \oplus \G_1^{\vee}(a) \to F(-c-m+a) \to 0.
$$
Then $h^0E(-c'+a)>0$ and a non trivial global section $s$ of $E(-c'+a)$ gives
a non zero map
$$
s:\Oc_X(-a) \to E(-c').
$$
Arguing as above this implies that $E$ splits which
is absurd.
\par\medskip
Then to finish the proof that case (i) can not arise we have to show the lemma.
\par\medskip
\noindent{\bf Proof of the Lemma.}
If $h^0 E(-a)>0$,
since
$E$ is normalized then $-a \geq 0$ i.e. $a \leq 0$.
Suppose that
$h^0F(-m-c+a)>0$. Since $F$ is normalized then
$-m-c+a \geq 0$. Then from conditions $-m-c+a \geq 0$ and $-a \geq 0$ we derive
condition
$c+m \leq 0$
which is absurd since by hypotheses we have condition $c+m>0$ (see the table).
\qed
\par\medskip
To show the theorem it remains now to consider the case (ii) i.e. when
$\G$ has an indecomposable
summand which is ACM of rank equal to $2$, i.e. when
\begin{equation} \label{eq:spl}
\G=\G_1 \oplus \G_2
\end{equation}
with both $\G_i$ ACM of rank equal to 2.
Of course, we may assume that $\G_i$ are
both indecomposable otherwise we reduce to the case (i) above.
Then we have non trivial extension class
\begin{equation} \label{eq:ext2}
\begin{CD}
0 @>{}>> F(m) @>f>> \G_2 \oplus \G_1 @>g>> E @>{}>> 0.
\end{CD}
\end{equation}

The extension class (\ref{eq:ext2}) is non trivial by assumption.
Moreover one has
\begin{equation} \label{eq:cl}
c_1(\G_i) \notin \{c_1(\F(m),c_1(E)\}
\end{equation}
for $i=1,2$,
by the corollary to Lemma 1.2.8 in \cite{OSS}. Indeed, suppose that
 $c_1(\G_i) \in \{c_1(F(m)),c_1(E)\}$ for
at least on $i=1,2$. Here we note that at least one of the bundles $F(m)$
and $E$ is stable. Hence $\G_i$'s are semistable and one of these is always stable.
Then from the above exact sequence we have map
between semistable bundles of the same rank with the same first Chern class
where at least one is stable. Therefore this map is an isomorphism and
hence the extension class is trivial, which is absurd.

Then to show that the splitting of (\ref{eq:spl}) can not arise we will use
Proposition \ref{prop:list} and a direct
computation on the Chern classes. It will show that the only possibility is that
the extension class (\ref{eq:ext1}) is trivial, which is absurd, since by assumption
$\G$ is represented by a non trivial class in $Ext^1(E,F(m))$.

To start with, we notice that
the bundle $\G$ of (\ref{eq:ext1}) is normalized since
so are $F$ and $E$, and $m \leq 0$. In particular also $\G_1$ and
$\G_2$ are normalized. \par
Then we consider all the possible splitting type of $\G$ under condition
(\ref{eq:cl}) in all the cases
(1)-(7) of the table.
Some of these decomposition are easy to show to be impossible, so we
give here only the cases, which require some more computations.
\par\smallskip
Case (1). In this case (see the table)
$$
(c_1(\G),c_2(\G))=(5,58).
$$
By Proposition \ref{prop:list} and by condition (\ref{eq:cl})
if $\G \cong \G_1\oplus \G_2$ splits then
$c_2(\G_1)=20$ and $c_2(\G_2)=\alpha$, with $\alpha=11,12,13,14$, are the
only possible cases.
A direct calculation on the Chern classes shows in these cases
$(c_1(\G),c_2(\G))\ne (5,58)$.\par
Case (2). In this case
$$
(c_1(\G),c_2(\G))=(2,18).
$$
By Proposition \ref{prop:list} and by condition (\ref{eq:cl}) if $\G \cong \G_1 \oplus \G_2$ then
we have only two possibilities: either
$(c_1(\G_1),c_2(\G_1))=(2,14)$ and $(c_1(\G_2),c_2(\G_2))=(0,4)$
or
$(c_1(\G_1),c_2(\G_1))=(2,13)$ and $(c_1(\G_2),c_2(\G_2))=(0,5)$.
The first case is impossible since by Riemann-Roch theorem
we have
$h^0F(-1)+h^0(E)=1<h^0(\G_1)+h^0(\G_2)=2$. The second case is also impossible
since one computes $h^0(\G_1)+h^0(\G_2)=3$.
\par
Case (3).
In this case we have
$$
(c_1(\G),c_2(\G))=(2,19).
$$
If $\G \cong \G_1 \oplus \G_2$ by Proposition
\ref{prop:list} and by condition (\ref{eq:cl}) we could have possible cases
$(c_1(\G_1),c_2(\G_1))=(1,6)$ and $(c_1(\G_2),c_2(\G_2))=(1,8)$ or
$(c_1(\G_1),c_2(\G_1))=(2,14)$ and $(c_1(\G_2),c_2(\G_2))=(0,5)$.
In the first case by Riemann-Roch we compute $0=h^0F(-1)+h^0E <
h^0\G_1+h^0\G_2=3$.
In the second case one concludes in similar way since $h^0\G_1+h^0\G_2=1$.
One concludes in similar way for the other cases.
\par
Case (4).
In this case we have
$$
(c_1(\G_1),c_2(\G_2))=(2,20).
$$
If $\G \cong \G_1 \oplus \G_2$ by Proposition
\ref{prop:list} and by condition (\ref{eq:cl}) we could have possible cases
$(c_1(\G_1),c_2(\G_1))=(1,\alpha)$ and $(c_1(\G_2),c_2(\G_2))=(1,\alpha')$
with $\alpha,\alpha'=4,6,8$ and $\alpha+\alpha'=15$ which is impossible.
\par
Cases (5)--(7). In this case we have
$$
(c_1(\G_1),c_2(\G_2))=(1,\alpha+8)
$$
where $\alpha=3,4,5$. If $\G \cong \G_1 \oplus \G_2$ by Proposition
\ref{prop:list} and by condition (\ref{eq:cl}) soon the conclusion follows.

\end{proof}


\begin{thebibliography}{ACGH}


\bibitem{AC} E.Arrondo and L.Costa,
\emph{Vector bundles on Fano 3-folds without intermediate cohomology},
Comm. Algebra \textbf{28} (2000), no. 8, 3899--3911.

\bibitem{BGS} R.O.Buchweitz, G.M.Greuel, and F.O.Schreyer,
\emph{Cohen-Macaulay modules on hypersurface singularities II},
Invent. Math.
\textbf{88} (1987), 165--182.

\bibitem{cm1} L.Chiantini and C.Madonna,
{\it ACM bundles on a general quintic threefold},
Matematiche (Catania) \textbf{55} (2000), no. 2, 239--258.

\bibitem{H1} R.Hartshorne,
\emph{Stable vector bundles of rank 2 on $\mathbb{P}^3$},
Math. Ann. \textbf{238} (1978), 229--280.

\bibitem{ks} S.Katz and S.Stromme,
\emph{Schubert}, a Maple package for intersection theory
and enumerative geometry, from
website http://www.mi.uib.no/schubert/

\bibitem{mad3} C.G.Madonna, {\it ACM bundles on prime Fano
threefolds and complete intersection Calabi Yau threefolds},
Rev. Roumaine Math. Pures Appl. \textbf{47} (2002), no.2, 211-222.

\bibitem{OSS} C.Okonek, M.Schneider and H.Spindler,
Vector bundles on complex projective spaces,
Progress in Mathematics \textbf{3}, 1980, pp.389.

\end{thebibliography}
\end{document}